\documentclass[a4paper,12pt]{article}

\usepackage{listings}
\usepackage{longtable}

\usepackage{graphicx}

\usepackage[font=scriptsize]{caption}

\usepackage{float}

\usepackage{amsmath,amssymb,mathrsfs,eucal}

\usepackage[T1,T2A]{fontenc}

\font\cyrfont=wncyss10

\def\sza{\hbox{\cyrfont X}} 

\newtheorem{thm}{Theorem}

\newtheorem{conj}[thm]{Conjecture}

\usepackage{color}

\begin{document}

\title{Critical $L$-values of Gross curves}

\author{Andrzej D\k{a}browski, Tomasz J\k{e}drzejak and Lucjan Szymaszkiewicz} 

\date{}

\maketitle{}

Let $K=\Bbb Q(\sqrt{-q})$, where $q$ is a prime congruent to $3$ modulo $4$. Let $A=A(q)$ 
denote the Gross curve \cite{Gr}. This $\mathbb Q$-curve is defined over $F=\mathbb Q(j)$ 
(where $j=j(\frac {1+\sqrt{-q}} {2}))$ and has good reduction outside $q$. 
$A(q)$ can be defined by the Weierstrass equation 
$$
y^2 = x^3 + {mq \over 2^33}x - {nq^2 \over 2^53^3}, 
$$
where $m$ and $n$ are the unique real numbers (they actually lie in $F$) such that 
$$
m^3=j, \; -n^2q=j-1728, \; sgn(n)=\left({2 \over q}\right). 
$$
Over $H=FK$ (the Hilbert class field of $K$) it acquires complex multiplication 
by $\mathbb {Z}[\frac {1+\sqrt{-q}} {2}]$.  It is known (\cite{Gr}, section 22) that the 
groups $A(F)$ and $A(H)$ are both finite if $q$ is a prime congruent to $7$ modulo $8$.

In this note we use Magma \cite{BCP} 
to calculate the values $L(E/H,1)$ for all such $q$'s up to some reasonable ranges 
for all primes $q$ congruent to $7$ modulo $8$.  
All these values are non-zero, and using the Birch and Swinnerton-Dyer conjecture, 
we can calculate hypothetical orders of $\sza(A/H)$ for all such $q$ up to $4831$. 
Our calculations agree with those given by Rodriguez Villegas \cite{RV}: he  
calculated conjectural orders of $\sza(A/F)$ by another method for all such $q$ up to $2927$, 
but one expects that $|\sza(A/H)| = |\sza(A/F)|^2$.

\bigskip

Now let us say a few words about the Magma implementation. 
The algorithm is analogous to those used by the authors in \cite{DJS}: it  uses the fact that 
 $L$-series of an elliptic curve over $H$ splits into factors 
corresponding to Grossencharacters twisted by Hilbert 
characters and its conjugates (so it uses the classical Hecke-Deuring theory 
linking elliptic curves with $CM$ to Grossencharacters, with keeping 
track of the effect of twisting).

\bigskip

Let us write down an explicit conjectural formula for the order of $\sza(E/H)$ 
(see \cite{Gr}, formula (21.3.1) on page 71).  Let $h$ denote the class number of $K$, 
and put $\epsilon_q(c):= \left(\frac c q\right)$.

\begin{conj} 
$$
\#(\sza(E/H)) = {q^{{3h \over 2}}  L(E/H,1) \over 
2^{3h-4} \pi^h  \prod_{0<c<q } 
\Gamma \left(\frac c q \right)^{\epsilon_q(c)}}.  
$$
\end{conj}

The above conjecture agrees with the Birch and Swinnerton-Dyer conjecture for $E/H$.

\bigskip 
 
John Coates informed one of us (A. D.) that one should be able to use the Iwasawa 
theory being developed in \cite{CKLT} to prove the above conjecture.

\bigskip  

Below we use Conjecture 1 to calculate $\#(\sza(E/H))$ (i.e. the analytic order of 
$\sza(E/H)$) for all primes $q$ congruent to $7$ modulo $8$ up to $4831$. 
Note that, if $L(E/H,1) \not=0$, then we know by Iwasawa theory that the Tate-Shafarevich 
group $\sza(E/H)$ is finite.

\begin{center}
\tiny
\begin{longtable}{|r|r|r|r|} 
\hline 
\multicolumn{1}{|c|}{$q$} 
& 
\multicolumn{1}{|c|}{$h$} 
&
\multicolumn{1}{|c|}{$\sqrt{L(E/H,1)}$} 
& 
\multicolumn{1}{|c|}{$\sqrt[4]{\#(\sza(E/H))}$} 
\\ 
\hline 
\endhead 
\hline 
\multicolumn{4}{|r|}{{Continued on next page}} \\
\hline
\endfoot

\hline 
\hline
\endlastfoot 

23 &	3 &	1.917306455589027597813347690641 &	$1$  \\  
31 &	3 &	1.023129835390262923123201119846 &	$1$  \\  
47 &	5 &	0.884007553190413369641732674593 &	$1$  \\  
71 &	7 &	2.480099873641678608433236547276 &	$3$  \\  
79 &	5 &	1.110556441918383299896561583115 &	$3$   \\  
103 &	5 &	0.038949744284901045024462605816 &	$1$   \\  
127 &	5 &	0.129399306674190803848276854132 &	$3$   \\  
151 &	7 &	0.342154370722521874022498924491 &	$3^2$   \\  
167 &	11 &	8.499196825809003222870613079977 &	$5 \cdot 23$   \\  
191 &	13 &	1.665085046661615396019042911840 &	$131$   \\  
199 &	9 &	0.676620807437586519563053934692 &	$43$   \\  
223 &	7 &	0.703326980964752529672024931405 &	$47$   \\  
239 &	15 &	0.453918696628391215678098573105 &	$3 \cdot 7 \cdot 17$  \\  
263 &	13 &	0.001257914062960834542480686368 &	$3 \cdot 7$  \\  
271 &	11 &	0.008805541008859093065049384872 &	$29$   \\  
311 &	19 &	0.648486521091286753635016640407 &	$3 \cdot 5 \cdot 479$   \\  
359 &	19 &	1.141338008694732473036144208271 &	$30529$   \\  
367 &	9 &	0.052953825151944195301873272368 &	$3^2 \cdot 19$   \\  
383 &	17 &	0.019680887700441358197318539018 &	$3^4 \cdot 41$   \\  
431 &	21 &	0.011456612474963829242166039245 &	$2^3 \cdot 3 \cdot 1249$   \\  
439 &	15 &	0.394270839749995230770685269658 &	$11311$   \\  
463 &	7 &	0.004174495009728004690701026903 &	$73$   \\  
479 &	25 &	0.048836267621336618889098778099 &	$353 \cdot 2273$   \\  
487 &	7 &	0.000000482116307460711439687440 &	$1$   \\  
503 &	21 &	0.000000020126315921323041547985 &	$151$   \\  
599 &	25 &	0.208905265188856310190917073528 &	$373 \cdot 52541$   \\  
607 &	13 &	0.003438572677705639141578890847 &	$3 \cdot 17 \cdot 73$   \\  
631 &	13 &	0.005900523952719352309760312826 &	$3^2 \cdot 911$   \\  
647 &	23 &	0.000045636135620621188056552483 &	$7 \cdot 33377$   \\  
719 &	31 &	0.046183960650516366354860538188 &	$1667632553$   \\  
727 &	13 &	0.001612218295199165920123560753 &	$3 \cdot 31 \cdot 101$   \\  
743 &	21 &	0.000091760752460157513132242927 &	$553649$   \\  
751 &	15 &	0.000814583878345215157582347921 &	$2^2 \cdot 7 \cdot 1087$   \\  
823 &	9 &	0.001813020189068166124586066929 &	$5 \cdot 577$   \\  
839 &	33 &	0.422592763639069040401160660086 &	$136363465939$   \\  
863 &	21 &	0.018883707188218542056256153867 &	$29 \cdot 1491911$   \\  
887 &	29 &	0.000483299398006109105446142215 &	$5 \cdot 11 \cdot 307 \cdot 47339$   \\  
911 &	31 &	0.212859814115063951127843593423 &	$2237 \cdot  45449137$  \\  
919 &	19 &	0.016284557379555919482008839655 &	$5 \cdot 37 \cdot 37747$   \\  
967 &	11 &	0.029126318364476542847711035736 &	$19 \cdot 5527$   \\  
983 &	27 &	0.007398840131454764434542833427 &	$4475472187$   \\  
991 &	17 &	0.000071072857379564759822178803 &	$3^2 \cdot 71 \cdot 509$   \\  
1031 &	35 &	0.015211980849109403793109862353 &	$37 \cdot 311 \cdot 6337 \cdot 25799$  \\  
1039 &	23 &	0.000355597434008835197149888602 &	$127 \cdot 499 \cdot 709$   \\  
1063 &	19 &	0.000237599218321558177783225700 &	$3 \cdot 277 \cdot 2999$   \\  
1087 &	9 &	0.000262116062494455330095732820 &	$59 \cdot 127$   \\  
1103 &	23 &	0.000000393293343069791012052711 &	$353 \cdot 36857$   \\  
1151 &	41 &	0.000003003732321349298811441864 &	$8422760498003$   \\  
1223 &	35 &	0.000002094475846449321713752880 &	$3 \cdot 5^2 \cdot 67 \cdot 101 \cdot 528413$   \\  
1231 &	27 &	0.000000489869434525845537217373 &	$229 \cdot 730753$   \\  
1279 &	23 &	0.000007115539976069114205127938 &	$5 \cdot 17 \cdot 29 \cdot 34549$   \\  
1303 &	11 &	0.000000671077509674452539881023 &	$3^2 \cdot 5 \cdot 7 \cdot 17$   \\  
1319 &	45 &	0.021858605815195499738879075748 &	$3^5 \cdot 124771 \cdot 4605889021$   \\  
1327 &	15 &	0.000000023399902732748609947789 &	$19 \cdot 809$   \\  
1367 &	25 &	0.000000045835215800037827192892 &	$2^4 \cdot 3 \cdot 5^2 \cdot 265883$   \\  
1399 &	27 &	0.000000177792482346179165660951 &	$2^2 \cdot 3 \cdot 17 \cdot 3006713$   \\  
1423 &	9 &	0.004136021545793286370439313380 &	$2^2 \cdot 19 \cdot 3163$   \\  
1439 &	39 &	0.000000002902257072083434474750 &	$5^3 \cdot 36575946427$   \\  
1447 &	23 &	0.000000004576434428595770264796 &	$3 \cdot 1464929$  \\  
1471 &	23 &	0.000009091947043380468173531523 &	$5 \cdot 4349 \cdot 20063$   \\  
1487 &	37 &	0.000000000054291128295627901658 &	$13 \cdot 1423 \cdot 1459 \cdot 6599$   \\  
1511 &	49 &	0.000000001232608907792342455037 &	$3^2 \cdot 109 \cdot 193 \cdot 2699 \cdot  18280687$  \\  
1543 &	19 &	0.000046300807064951614629771602 &	$5 \cdot 12952607$   \\  
1559 &	51 &	0.000000001909704251728812403863 &	$67075913 \cdot 1541369429$   \\  
1567 &	15 &	0.000269619672929093678586912043 &	$3 \cdot 1471 \cdot 2803$   \\  
1583 &	33 &	0.000591748987259341323505132764 &	$7369 \cdot 12135601151$   \\  
1607 &	27 &	0.000051925063264633860691291385 &	$524286640211$   \\  
1663 &	17 &	0.000001469043900569342928712801 &	$3 \cdot 13 \cdot 121609$   \\  
1759 &	27 &	0.002034329386922218553153432265 &	$5 \cdot 1697 \cdot 9883 \cdot 15881$   \\  
1783 &	17 &	0.000000397460901878119350764458 &	$3^3 \cdot 7 \cdot 19 \cdot 2267$   \\  
1823 &	45 &	0.011540308508217040309927578478 &	$2^8 \cdot 823 \cdot  675413 \cdot 349590427$  \\  
1831 &	19 &	0.000001902416570890513127395632 &	$5 \cdot 23 \cdot 31 \cdot 139 \cdot 431$   \\  
1847 &	43 &	0.000000006214811080409277421457 &	$23 \cdot 167 \cdot  486769 \cdot  6302473$  \\  
1871 &	45 &	0.035640901248044787979620375731 &	$3^3 \cdot 7 \cdot 1740181568221388179$  \\  
1879 &	27 &	0.000000710833676781376253894792 &	$2^2 \cdot 3^2 \cdot 17 \cdot 172144933$   \\  
1951 &	33 &	0.000003529625675343434758248008 &	$65599 \cdot  546118759$  \\  
1999 &	27 &	0.000000005979473199318522728872 &	$53 \cdot 321482179$   \\  
2039 &	45 &	0.000000002182100425038097541281 &	$67 \cdot 528667 \cdot  17305594861$  \\  
2063 &	45 &	0.000038203268549928881568814494 &	$2^4 \cdot 5 \cdot 73 \cdot 89^2 \cdot 1004027 \cdot 1152643$  \\  
2087 &	35 &	0.000000000000014489648203923018 &	$17 \cdot 113 \cdot 241 \cdot 832621$   \\  
2111 &	49 &	0.000020147662461750776760561609 &	$5^2 \cdot 131 \cdot 9533 \cdot 103725964580429$   \\  
2143 &	13 &	0.000000073692751127076619011376 &	$3 \cdot 5 \cdot 54773$   \\  
2207 &	39 &	0.000000034190476384724778836118 &	$42907145154487223$   \\  
2239 &	35 &	0.000000000114240665575503380052 &	$19^2 \cdot  40105774091$  \\  
2287 &	29 &	0.000000019569362107701266283239 &	$526738119443$   \\  
2311 &	29 &	0.000017028956170501832053023258 &	$1390343 \cdot  28260481$  \\  
2351 &	63 &	0.000000000004962696845938121650 &	$11 \cdot 1171 \cdot 3739 \cdot 85780523 \cdot  788963191$  \\  
2383 &	29 &	0.000000000018701167333899812278 &	$7 \cdot 11 \cdot 29 \cdot  14437343$  \\  
2399 &	59 &	0.000000067762923493270317325060 &	$167 \cdot 397 \cdot  100609 \cdot 4027834961600243$  \\  
2423 &	33 &	0.000000000004118523178705156894 &	$199 \cdot  139095471773$  \\  
2447 &	37 &	0.000000000003265549315360713458 &	$3 \cdot 7^2 \cdot 17 \cdot 23 \cdot 82067 \cdot 165293$   \\  
2503 &	21 &	0.000000013421414398265685521548 &	$3^3 \cdot 233 \cdot 219169$   \\  
2543 &	35 &	0.000000000061443201800583958759 &	$2^3 \cdot 3 \cdot 5 \cdot 47 \cdot 1877 \cdot 167286883$   \\  
2551 &	41 &	0.000000516944886355770251892511 &	$2143 \cdot 777478894428461$   \\  
2591 &	57 &	0.000000000000041168041867756552 &	$31 \cdot 11617 \cdot  74463344072813299$ \\  
2647 &	15 &	0.000000104412115837007034423836 &	$3 \cdot 53 \cdot 149 \cdot 2819$   \\  
2663 &	43 &	0.000000038668404124883865748112 &	$13 \cdot 514531 \cdot 24237844896113$  \\  
2671 &	23 &	0.000000000000006590944134402918 &	$3 \cdot 23894833$   \\  
2687 &	51 &	0.000000000063446002820259446740 &	$2^2 \cdot 4813093 \cdot  324237065387803$  \\  
2711 &	53 &	0.000000000000000102784949519600 &	$43^2 \cdot 109 \cdot 759144336107881$   \\  
2719 &	41 &	0.000000000112187841555263134829 &	$5 \cdot 7 \cdot 83059 \cdot 31620028637$   \\  
2767 &	21 &	0.000000000061044792900945221648 &	$2^2 \cdot 3 \cdot 41^2 \cdot 43 \cdot 461$   \\  
2791 &	39 &	0.000000000005517487082621520322 &	$607 \cdot  13252797441179$  \\  
2879 &	57 &	0.000000000003933229323617751788 &	$11 \cdot 199 \cdot  3671245687 \cdot  777206990563$ \\  
2887 &	25 &	0.000014317579969749210837178220 &	$821 \cdot 106357 \cdot 266801$   \\  
2903 &	59 &	0.000000397671660764977706938361 &	$3^2 \cdot 5 \cdot 7 \cdot 61 \cdot 82891 \cdot 77113721 \cdot   36485503729$  \\  
2927 &	31 &	0.000000000053108304414288284399 &	$41 \cdot 53 \cdot 1723 \cdot  307885729$  \\  
2999 &	73 &	0.000000000000190277589599147568 &	$3 \cdot 5 \cdot 563 \cdot  648723344089 \cdot 2021222121768757$  \\  
3023 &	47 &	0.000000000008673013112970432807 &	$33070861 \cdot  59499888244757$  \\  
3079 &	41 &	0.000000000054487305510806875829 &	$372451 \cdot 3108015712297$   \\  
3119 &	69 &	0.014845943183323009021216005907 &	$2^2 \cdot 3 \cdot 11 \cdot 271 \cdot 2113759 \cdot 277103227 \cdot 7719400992115643$  \\  
3167 &	53 &	0.000000000419958531119441000457 &	$3 \cdot 7^2 \cdot 17 \cdot 202059233 \cdot 22795690268033$   \\  
3191 &	69 &	0.000000073648283059495109159921 &	$9460962968276441 \cdot 160726061974992413$  \\  
3271 &	27 &	0.000000000004820011697320255318 &	$2^2 \cdot 3^3 \cdot 1321 \cdot 1559 \cdot 8521$  \\  
3319 &	41 &	0.000000000000000043650734177687 &	$3 \cdot 5 \cdot 11 \cdot 14525383286597$  \\  
3343 &	19 &	0.000000000465261296224572289673 &	$3^2 \cdot 5^2 \cdot 14609197$  \\  
3359 &	69 &	0.000000002215874222576132135412 &	$109 \cdot 14821 \cdot 47207 \cdot 562021 \cdot 5120491 \cdot 8211851209$   \\  
3391 &	37 &	0.000000012667841601214216959683 &	$3^3 \cdot 29 \cdot 251 \cdot  20380307783611$ \\  
3407 &	57 &	0.000000000000069504074168580772 &	$2^2 \cdot 83 \cdot 6607 \cdot 1805299 \cdot  13807764284591$  \\  
3463 &	19 &	0.000000000001451107621933911886 &	$5 \cdot 7^3 \cdot 13 \cdot 59 \cdot 571$  \\  
3511 &	41 &	0.000000000253148503966522826850 &	$5 \cdot 4502030527785827207$  \\  
3527 &	65 &	0.000000000007634089265859983699 &	$3 \cdot 43 \cdot 397 \cdot  115891 \cdot  141941 \cdot 152539  \cdot 17046231029$  \\  
3559 &	45 &	0.000000000000000006597486533662 &	$2^2 \cdot 41 \cdot 179 \cdot 521 \cdot 9377 \cdot 3906323$  \\  
3583 &	29 &	0.000000000001148486625059495068 &	$19 \cdot 653277444511$  \\  
3607 &	19 &	0.000000001553716006672287458844 &	$83663 \cdot  223099$ \\  
3623 &	45 &	0.000000000000919701344137217043 &	$3 \cdot 59 \cdot 40648887911215217617$  \\  
3631 &	43 &	0.000000000022103794221839759633 &	$19 \cdot 337 \cdot 397 \cdot 13339 \cdot 6057046621$  \\  
3671 &	81 &	0.000000005729179615089549879246 &	$5 \cdot 17 \cdot 23 \cdot 93637 \cdot 404309 \cdot 7413852725693 \cdot 10527977612261$  \\  
3719 &	67 &	0.000000003987873570890300692631 &	$4621 \cdot 9879817 \cdot  32792616424035435083$ \\  
3727 &	31 &	0.000000000000000715924299282265 &	$3 \cdot 521 \cdot 2234679487$  \\  
3767 &	39 &	0.000000000000000000023229393826 &	$3 \cdot 11 \cdot 47 \cdot 59 \cdot  159124023689$ \\  
3823 &	29 &	0.000000000011639491490795636800 &	$3 \cdot 5 \cdot 61 \cdot 647 \cdot 75001681$ \\  
3847 &	23 &	0.000000033915083714883628181974 &	$3^3 \cdot 3727 \cdot 111765767$  \\  
3863 &	61 &	0.000000266167237138593593160127 &	$3 \cdot 5206955446969 \cdot 14572866296207015941$  \\  
3911 &	83 &	0.000000014026513503239094034158 &	$89381 \cdot 10120657 \cdot 26771077 \cdot 34445368269110070820019$  \\  
3919 &	39 &	0.000000000000444262371239767809 &	$2^2 \cdot 5^2 \cdot  45310663031163427$ \\  
3943 &	27 &	0.000000025420324493979393463674 &	$173 \cdot  2349245993371$  \\  
3967 &	33 &	0.000000105745446379198858088646 &	$2^2 \cdot 2287 \cdot  114728185930409$ \\  
4007 &	57 &	0.000000000000201960850286850752 &	$3^2 \cdot 137 \cdot 227 \cdot 331 \cdot 599  \cdot 255045098344833101$ \\  

4079 &	85 &	0.000000000000000082140602546709 &	$2^4 \cdot 3 \cdot 7 \cdot 11 \cdot 614594868178532601811498841258583389$ \\ 
4111 &	39 &	0.000000000000264993482744019345 &	$3^2 \cdot 5 \cdot 97 \cdot 1167954319891727$ \\ 
4127 &	49 &	0.000000000000193127996090472079 &	$3 \cdot 4591 \cdot 58067 \cdot 8556016692999833$ \\ 
4159 &	31 &	0.000000002334840823717447231766 &	$7 \cdot 5419 \cdot 6698995552091$ \\ 
4231 &	51 &	0.000000000002854189047798681349 &	$503 \cdot 1933 \cdot 9687106562299193749$ \\ 
4271 &	65 &	0.000000000000022249456103412573 &	$3 \cdot 257 \cdot 416400721169 \cdot 
1259073376623055453$ \\  
4327 &	19 &	0.000361493283256042320918368116 &	$407437 \cdot 294345721$ \\ 
4391 &	79 &	0.000000000000019695220419682544 &	$3^2 \cdot 5 \cdot 13 \cdot 173 \cdot 10592464484160289704717797797701109$ \\ 
4423 &	33 &	0.000000001925109169211800528055 &	$2^2 \cdot 378508129 \cdot 823095127$ \\ 
4447 &	17 &	0.000000000001867493388178899269 &	$29 \cdot 79 \cdot 1041449$ \\ 
4463 &	55 &	0.000000000000984566972390265807 &	$3 \cdot 47 \cdot 485739683983935166227585053$ \\ 
4519 &	29 &	0.000000000000000453913600988310 &	$3 \cdot 13 \cdot 19 \cdot 107166297101$ \\ 
4567 &	33 &	0.000000000003742086161710712616 &	$5 \cdot 919 \cdot 40002808124563$ \\ 
4583 &	61 &	0.000000000110718923407311375654 &	$337 \cdot 8320811 \cdot 727967861399585028840511$ \\ 
4591 &	49 &	0.000000000000164991042996831997 &	$5 \cdot 451733284465625717690381$ \\ 
4639 &	51 &	0.000000000000000262582321262483 &	$2^8 \cdot 3 \cdot 1900987826216667061231$ \\ 
4663 &	33 &	0.000000000000003556164981697977 &	$2^4 \cdot 5 \cdot 7 \cdot 3015216230623$ \\  

4679 &          91 &   0.000000000025511112078436601931 &        $2^6 \cdot 3 \cdot  36389  \cdot  30546984990724085865497148667129862189333$ \\   

4703 &         75  &  0.000000000003892521768891042131  &        $2^2 \cdot 37 \cdot 61^2 \cdot 69380888539 \cdot  24421524889725136602727$ \\  

4751 &         91  &  0.000000821131470506302985971524 &        $53 \cdot 263 \cdot   107033  \cdot 51896369002543288357  \cdot  1315458720324753637403$  \\  

4759 &         55  &  0.000000000000004952618783154918 &        $13693553 \cdot 63045641606597183767$  \\  

4783 &         23  &  0.000000004277461191967938167424 &        $3 \cdot 7 \cdot 13 \cdot 23^2 \cdot 43 \cdot 181 \cdot 115891$
 \\ 
 
4799 &         63  &  0.000000000000000000063775415526 &        $2^2 \cdot 3^5 \cdot 1093 \cdot 7955015043561253971333769$  \\  

4831 &         33  &  0.000000000000003424736207358667 &        $29 \cdot 419 \cdot 6691 \cdot 988157509$  \\

\end{longtable}
\end{center}

We  thank John Coates for suggesting the problem to us, and for very inspiring correspondence.

Institute of Mathematics, University of Szczecin, Wielkopolska 15, 
70-451 Szczecin, Poland; E-mail addresses: andrzej.dabrowski@usz.edu.pl and dabrowskiandrzej7@gmail.com;  
tjedrzejak@gmail.com;   lucjansz@gmail.com

\end{document}